\documentclass[reqno,11pt]{amsart}

\newtheorem{theorem}{Theorem}[section]
\newtheorem{lemma}[theorem]{Lemma}
\newtheorem{corollary}[theorem]{Corollary}

\theoremstyle{definition}
\newtheorem{assumption}[theorem]{Assumption}

\theoremstyle{remark}
\newtheorem{remark}[theorem]{Remark}

\makeatletter
\def\dashint{\operatorname%
{\,\,\text{\bf--}\kern-.98em\DOTSI\intop\ilimits@\!\!}}
\makeatother

\newcommand\bR{\mathbb{R}}

\newcommand\cF{\mathcal{F}}

\newcommand{\osc}{{\rm osc}\,}
\newcommand{\loc}{{\rm loc}\,}

\newcommand{\WO}{\overset{\scriptscriptstyle0}%
{W}\,\!}

 \newcommand{\mysection}[1]{\section{#1}
 \setcounter{equation}{0}}

\newcommand{\nliminf}{\operatornamewithlimits{\underline{lim}}}

\begin{document}

\title%
{Elliptic equations with VMO a,
  b$\,\in L_{d}$, and c$\,\in L_{d/2}$}
\author{N.V. Krylov}

\email{nkrylov@umn.edu}
\address{127 Vincent Hall, University of Minnesota, Minneapolis, MN, 55455}
 
 

\begin{abstract} We consider elliptic equations with operators
$L=a^{ij}D_{ij}+b^{i}D_{i}-c$ with $a$ being almost in VMO,
$b\in L_{d}$ and $c\in L_{q}$, $c\geq0$, $d>q\geq d/2$. We prove
the solvability of $Lu=f\in L_{p}$ in bounded $C^{1,1}$-domains,
 $1<p\leq q$,
and of $\lambda u-Lu=f$ in the whole space for any 
$\lambda>0$.
 Weak uniqueness of the martingale
problem associated with such operators is also obtained. 
\end{abstract}

\maketitle

\mysection{Introduction}
                                                  \label{section 3.11.1}
Let $\bR^{d}$ be a $d-$dimensional Euclidean space of points
$x=(x^{1},...,x^{d})$ with $d\geq2$.
We are dealing with a uniformly  elliptic operator
  $$
Lu( x)= a^{ij}( x)D_{ij}u ( x)+
b^{i}( x)D_{i}u( x)-c(x)u(x) ,\quad D_{i}=\frac{\partial}{\partial x^{i}},
\quad D_{ij}=D_{i}D_{j},
$$
acting on functions given on $\bR^{d}$. Throughout the article the
numbers
$p,q\in(1,\infty)$ are fixed and assumed to satisfy 
either
\begin{equation}
                                                \label{3.9.1}
d/2<q<d,\quad 1<p\leq q,
\end{equation}
or
\begin{equation}
                                                \label{3.9.2}
q=d/2,\quad 1<p< d/2\quad (\text{and}\quad d\geq3).
\end{equation}

We assume that $b\in L_{d}(\bR^{d})$ and $c\in L_{q}(\bR^{d})$, $c\geq 0$.
Note that the case that $q=d/2$ is   generally  not excluded. 
 However, it is excluded if $d=2$ because there are no $p$
satisfying $1<p<1$. If $d=2$ we should have $c\in L_{q}(\bR^{d})$
with $q>1$. 
We also assume that $a$ is bounded and almost in VMO and prove the unique
solvability results for the equation $Lu=f\in L_{p}(G)$ in regular domains $G$
in the class $\WO^{2}_{p}(G)$ and for the equation
$\lambda u-Lu=f$ in $\bR^{d}$ for any $\lambda>0$ in the class
$W^{2}_{p}(\bR^{d})$. We apply these results to prove that
the corresponding solutions of It\^o's stochastic equation possess the weak
uniqueness property.

To the best of the author's knowledge these results are new
even if $a^{ij}=\delta^{ij}$, however, much work was done
in this case.

 G. Stampacchia in \cite{St_65} (1965) was probably the first
author who presented the $W^{1}_{2}$-solvability
 theory of divergence form equations
with $b\in L_{d}(\bR^{d})$ and $c\in L_{d/2}(\bR^{ d })$,
with some additional restrictions on $c$ but
without assuming the smallness of the norms of $b$ and $c$ as well as
  without assuming that the domain $G$ in which the equation
is solved is small. The restriction on $c$ can be summarized
as follows (see Theorem 3.4 in \cite{St_65}): $c=c'+\lambda$,
where the parameter $\lambda$ should not belong to a countable set,
which is known to be lying below some $\bar \lambda>0$.
There is a plethora of other important results in \cite{St_65},
but we will discuss only the one mentioned above which is most related
to our own results in case $a^{ij}=\delta^{ij}$. The free term 
\cite{St_65} is taken
in the divergence form
$f=D_{i}g_{i}$, where $g_{i}\in L_{2}$. To match this with
our $f\in L_{p}$ under condition \eqref{3.9.2} we have to have
$p>2d/(d+2)$ (and $d\geq 3$). Then in this level-ground situation
we have a solution $u\in W^{2}_{p}$ and   \cite{St_65}
guarantees only $u\in W^{1}_{2}$. At the same time $W^{2}_{p}
\subset W^{1}_{r}\subset W^{1}_{2}$ for an $r>2$. By the way,
in the statement of Theorem 3.4 of \cite{St_65} the condition that $d\geq 3$
is not included, but it is actually tacitly imposed (see pages 200-201
there where embedding theorems are applied).

The estimates leading to Theorem 3.4 of \cite{St_65}
are also found in 
O.A. Ladyzhenskaya and N.N. Ural'tseva book \cite{LU_73} (1973),
see pages 189-191 there, where the condition $d\geq 3$
is explicitly imposed. 

N.S. Trudinger in \cite{Tr_73} (1973) in the setting of
generally degenerate divergence
form operators, among many other things, removed the condition
on $c$ in \cite{St_65} and replaced it with just $c\geq0$.

In a recent article \cite{KK_19}
 Byungsoo Kang and Hyunseok Kim present a deep investigation
of divergence type equations with $b\in L_{d}$ (no $c$) and solutions
in $W^{1}_{p}$ and, if $a$ is more regular, in $W^{2}_{p}$. 

It is also worth mentioning the recent article \cite{ANPS_19} 
by  Apushkinskaya,   Nazarov,   Palagachev, and
Softova, in which
  nondivergence  
form equations with VMO $a$, $b\in L_{d}$, and $c\in L_{d/2}$
are considered and a priori estimates are obtained but on the right
 in these estimates the zeroth order norm of the unknown
function is present. 

Saisai Yang and Tusheng Zhang in \cite{YZ_18} (2018)
  use   probabilistic approach to prove, among other things,
 that there exists
a unique
$C^{1+\alpha}$-solution to the Dirichlet boundary 
 value problem   in case $a^{ij}=\delta^{ij}$
and $b$ and $c$ are signed measures
of Kato class $K_{d,1-\alpha}$, $\alpha\in(0,1)$. Although the pointwise
regularity of solutions in \cite{YZ_18} is stronger than ours $u\in C^{2-d/q}$
(see Corollary \ref{corollary 3.5.1}), it is worth mentioning that
a general $f\in L_{r}$ is in $K_{d,1-\alpha}$ only if $(1-\alpha) r>d$.
Therefore, generally $b\in L_{d}$ is not in any $K_{d,1-\alpha}$
and $c\in L_{d/2}$ is way out of $K_{d,1-\alpha}$. Therefore,
 the results of \cite{YZ_18} are not applicable in our case.

 The above discussion seems to support that even in the case 
of $a^{ij}=\delta^{ij}$ our PDE results were unknown. We prove them when $a\in
BMO$. Regarding numerous issues for equations with BMO main coefficients
and bounded lower order coefficients we refer the reader to
\cite{DK_11} and the references therein.

In what concerns the weak uniqueness os solutions of the
corresponding stochastic equation with drift in Lebesgue
spaces, much work has been done
mostly in the time nonhomogeneous case $b=b(t,x)$, mostly when
$a^{ij}=\delta^{ij}$. A good source of recent results and bibliography is the
paper by
 L. Beck, F. Flandoli, M. Gubinelli, and M. Maurelli
\cite{BFGM_19} (2019). In this paper the authors prove (see there
Theorem 5.4) an existence and uniqueness (stronger than the pathwise
uniqueness) theorem
applicable to our situation ($a^{ij}=\delta^{ij}$, $b\in L_{d}$) but only for
solutions with initial starting point having density and only 
in a class of solutions possessing certain properties. We prove
{\em weak\/} uniqueness but for any solution starting from any fixed point
and our $a\in BMO$. From the probabilistic point of view this is
also a new result complementing  the information in \cite{Kr04}.

\mysection{Equations in bounded domains}
                                                  \label{section 3.11.2}

 Fix numbers $\delta\in(0,1)$
and $\|b\|,\|c\| \in[0,\infty)$.

\begin{assumption}
                                                  \label{assumption 3.1.1}

The coefficients of $L$
   are measurable, the matrices $a(x)=(a^{ij}(x))$ are
symmetric and satisfy
\begin{equation}     
                            \label{7.18.1}
\delta^{-1}|\lambda|^{2}\geq
a^{ij}(x) \lambda^{i}\lambda^{j}\geq\delta|\lambda|^{2}
\end{equation}
for all $\lambda,x\in\bR^{d}$.
Also $c\geq 0$,
 $$
\|b\|_{ L_{d}(\bR^{d})}\leq\|b\| ,\quad 
\|c\|_{ L_{q}(\bR^{d})}\leq\|c\|.
$$
\end{assumption}

To state one more assumption we set $B_{r}(x)$ to be the open ball 
in $\bR^{d}$ of radius $r$ centered at $x$, $B_{r}=B_{r}(0)$.
Denote
$$
\osc (a,B_{\rho}(x))=
 |B_{\rho}|^{-2}
 \int_{y,z\in B_{\rho}(x)}|a( y)-a( z)|\,dydz,
$$
$$
a^{\# }_{r}=\sup_{ x \in\bR^{d }}\sup_{\rho< r}
\osc (a,B_{\rho}(x)) .
$$
 
 Set
$$
L_{0}u=a^{ij} D_{ij}u .
$$

Fix a bounded domain $G\subset\bR^{d}$ of class $C^{1,1}$.
Here is a particular case of Theorem 8 of \cite{DK_11},
 in which $Du$ is the gradient of $u$ and $D^{2}u$
is its Hessian. 
\begin{lemma}
                          \label{lemma 2.19.2}
Under Assumption  \ref{assumption 3.1.1}
for any $s\in(1,\infty)$   there exists $\theta_{0}=\theta_{0}(d,\delta,s)$
such that, if there is $r_{0}>0$ for which $a^{\#}_{r_{0}}\leq\theta_{0}$,
then there exist
  $\lambda_{0}\geq 1, N_{0}$,
depending only on $d,\delta,s,r_{0}$, and $G$,
  such that, for any $u\in \WO^{2}_{s}(G)$
and $\lambda\geq \lambda_{0}$,
\begin{equation}
                                 \label{2.19.3}
\|D^{2}u\|_{L_{s}( G )}+\sqrt\lambda\| Du\|_{L_{s}(G)}+\lambda\|
u\|_{L_{s}(G)}
\leq N_{0} \|L_{0}u-\lambda u\|_{L_{s}(G)}.
\end{equation}
 Furthermore, for any $f\in L_{s}(G)$
there exists a unique $u\in \WO^{2}_{s}(G)$
such that $L_{0}u-\lambda u=f$.
\end{lemma}

We fix $r_{0}>0$ and impose the following.

\begin{assumption}[$ p,r_{0}$] 
                                               \label{assumption 2.20.1}  
We have 
 $a^{\#}_{r_{0}}\leq\theta_{0}(d,\delta,p)$, where $\theta_{0}$
is taken from Lemma \ref{lemma 2.19.2}.
\end{assumption}

  Recall that we write $a\in VMO$ if 
$a^{\# }_{r}\to0$ as $r\downarrow 0$. So, our $a$ is ``almost'' in VMO.
It is also worth mentioning that $a\in VM0$ if, for instance,
($a$ is bounded and) $Da\in L_{d}(\bR^{d})$. An example of such
(uniformly nondegenerate bounded) $a$
is given by $2+I_{x \ne0}\zeta(x )\sin(\ln|\ln |x |)$,
where $\zeta$ is any smooth symmetric $d\times d$-matrix valued 
function vanishing for $|x|>1/2$ and
satisfying $|\zeta|\leq 1$.

Below by $\lambda_{0}$ we mean the one from Lemma
\ref{lemma 2.19.2} for $s=p$.

\begin{theorem}
                                                    \label{theorem 2.20.1}

Under Assumptions \ref{assumption 3.1.1} and \ref{assumption 2.20.1}
($p,r_{0}$)
  introduce the constant $N^{*}= N^{*}(p,q,d,G)$ as the best
constant such that
$$
\|Du\|_{L_{pd/(d-p)}(G)}+\|u\|_{L_{pq/(q-p)}(G)}
\leq N^{*}( \|D^{2}u\|_{L_{p}(G)} +\|  u\|_{L_{p}(G)})
$$
for any $u\in \WO^{2}_{p}(G)$. 
Assume that
\begin{equation}     
                            \label{2.20.3}
2N_{0}N^{*}(\|b\|_{L_{d}(G)}+\|c\|_{L_{q}(G)})\leq1 .
\end{equation}

Then for any $u\in \WO^{2}_{p}(G)$
and $\lambda\geq \lambda_{0}$,
\begin{equation}
                                 \label{2.20.4}
\|D^{2}u\|_{L_{p}(G)}+\sqrt\lambda\| Du\|_{L_{p}(G)}+\lambda\| u\|_{L_{p}(G)}
\leq 2N_{0} \|L u-\lambda u\|_{L_{p}(G)}.
\end{equation}
 Furthermore, for any $f\in L_{p}(G)$
there exists a unique $u\in \WO^{2}_{p}(G)$
such that $L u-\lambda u=f$.
\end{theorem}

Proof.  To prove \eqref{2.20.4} observe that
$$
N_{0} \|L _{0}u-\lambda u\|_{L_{p}( G )}
\leq N_{0} \|L u-\lambda u\|_{L_{p}( G )}+N_{0}(\|\,|b|\,|Du|\,\|
_{L_{p}( G )}+\|cu\|_{L_{p}( G )}),
$$
where the last term by H\"older's inequality,
embedding theorems, and \eqref{2.20.3} is less than
$$
N_{0}\|b\|_{L_{d}( G )}\|Du\|_{L_{pd/(d-p)}( G )}
+N_{0}\|c\|_{L_{q}( G )}\| u\|_{L_{pq/(q-p)}( G )}
$$
$$
\leq
N_{0}(\|b\|_{L_{d}( G )}+\|c\|_{L_{q}( G )})
N^{*}(\|D^{2}u\|_{L_{p}( G )}+\|  u\|_{L_{p}(G)})
$$
$$
\leq(1/2)(\|D^{2}u\|_{L_{p}( G )}+\|  u\|_{L_{p}(G)}).
$$
This shows that \eqref{2.20.4} follows from \eqref{2.19.3}.
The existence assertion of the theorem follows
as usual by the method of
continuity. The theorem is proved.

\begin{remark}
                                                  \label{remark 2.20.2}
The above estimates show that the operator $L$ is bounded
as an operator from $W^{2}_{p}( G )$ to $L_{p}( G )$ 
as long as $b\in L_{d}( G )$ and $c\in L_{q}(G)$.
\end{remark}

Next for our fixed $b$ and $c$ there exist  $b_{0},c_{0}\geq 0$ such that  
\begin{equation}     
                            \label{2.20.5}
2N_{0}N^{*}(\|bI_{|b|\geq b_{0}}\|_{L_{d}( G )}
+\|cI_{c\geq c_{0}}\|_{L_{q}( G )})\leq1 .
\end{equation}

\begin{theorem}
                                     \label{theorem 10.17.1}
Under Assumptions \ref{assumption 3.1.1} and \ref{assumption 2.20.1}
($p,r_{0}$)
there exist
  $\lambda_{1}\geq1, N $,
depending only on $d, \delta,p , r_{0}$, $b_{0}$, $c_{0}$, and $G$,
  such that, for any $u\in \WO^{2}_{p}(G)$
and $\lambda\geq \lambda_{1}$,
\begin{equation}
                                 \label{3.1.1}
\|D^{2}u\|_{L_{p}(G)}+\sqrt\lambda\| Du\|_{L_{p}(G)}+\lambda\|
u\|_{L_{p}(G)}
\leq N  \|L u-\lambda u\|_{L_{p}(G)}.
\end{equation}
 Furthermore, for any $f\in L_{p}(G)$
there exists a unique $u\in \WO^{2}_{p}(G)$
such that $L u-\lambda u=f$.
\end{theorem}

Proof. As usual, it suffices to prove the a priori
estimate \eqref{3.1.1}. By Theorem \ref{theorem 2.20.1}
its left hand side is dominated by
$$
2N_{0}\|L u-\lambda u\|_{L_{p}(G)}+2N_{0}
\|b^{i}D_{i}u I_{|b|\leq b_{0}}+cuI_{c\leq c_{0}}\|_{L_{p}(G)},
$$
where the last term, by interpolation
inequalities is less than
$$
N(\| D u  \|_{L_{p}(G)}+\|   u  \|_{L_{p}(G)})\leq (1/2)\|D^{2}u\|_{L_{p}( G )}
+N_{1}\| u\|_{L_{p}( G )}.
$$
This yields \eqref{3.1.1} for $\lambda\geq 2N_{1}$
and proves the theorem.

We denote the solution from Theorem \ref{theorem 10.17.1}
by $R_{\lambda+c}f$.
 
\begin{remark}
                                                  \label{remark 2.21.1}

By taking here $\lambda=\lambda_{1}$ in \eqref{3.1.1}
we see that for the same kind of $N$ as in \eqref{3.1.1}
and any $u\in \WO^{2}_{p}( G )$
\begin{equation}
                                               \label{10.22.2}
\| u\|_{W^{2}_{p}( G )}\leq N\big(\|Lu\|_{L_{p}( G )}
+\|u\|_{L_{p}( G )} \big).
\end{equation}
\end{remark}

 The next result, the proof of which is left to the reader,
is a standard consequence of Theorem \ref{theorem 10.17.1}
\begin{theorem}
                                                \label{theorem 3.10.1}
 Let   $a^{n},b^{n},c^{n}$, $n= 1,2,...$,
be a sequence of symmetric $d\times d$-matrix valued,
$\bR^{d}$-valued, and $[0,\infty)$-valued, respectively,
 measurable functions, satisfying Assumptions \ref{assumption 3.1.1}
 and \ref{assumption 2.20.1} ($p,r_{0}$)
(with the same $\delta$, $\|b\|,\|c\|$, and $\theta_{0}$  as above).
Let $f,f^{n}\in L_{p}(G)$ and suppose that
$ a^{n}\to a$ on $\bR^{d}$ (a.e.) and
$$
\|b-b^{n}\|_{L_{d}(G)}+\|c^{n}-c\|_{L_{q}(G)}
+\|f^{n}-f\|_{L_{p}(G)}\to0
$$
 as $n\to\infty$.
Let $\lambda\geq\lambda_{1}$, where $\lambda_{1}$
is taken from Theorem \ref{theorem 10.17.1}, and introduce $u^{n}$
as  unique $\WO^{2}_{p}(G)$-solutions of $\lambda u^{n}-L^{n}u^{n}=f $,
where the operator $L^{n}$ is constructed from  $a^{n},b^{n},c^{n}$.
Then 
$$
\lim_{n\to\infty}\|u^{n}-R_{\lambda+c}f\|_{W^{2}_{p}(G)}=0.
$$
\end{theorem}

By using mollifiers and properties of solutions of equations
with smooth coefficients we easily arrive at the following.
\begin{corollary}
                                         \label{corollary 3.10.1}
Under Assumptions \ref{assumption 3.1.1} and \ref{assumption 2.20.1}
 ($p,r_{0}$)
for $\lambda\geq\lambda_{1}$, where $\lambda_{1}$
is taken from Theorem \ref{theorem 10.17.1}, and any $f\in L_{p}(G)$
we have $|R_{\lambda+c}f|\leq R_{\lambda+c}|f|\leq R_{\lambda}|f|$
(a.e.).

\end{corollary}

Next we turn to some  properties of equations 
with $b\in L_{d}$ and $c\in L_{q}$. The main goal of these further
results, important in their own rights,
 is to prepare the necessary tools to be
able to treat the equations in the whole space for  any  $\lambda>0$
and in domains when $\lambda=0$.
 
\begin{lemma}
                                               \label{lemma 2.25.1}
Under Assumption \ref{assumption 3.1.1}
let $0<R\leq R_{0}<\infty$,
$\varepsilon\in(0,1]$, 
\begin{equation}
                                                           \label{3.9.3}
 d\geq t\geq p,\quad \frac{d}{p} <  1+\frac{d}{t},
\end{equation}
$u\in W^{2}_{p}( G )$,
$\zeta\in C^{\infty}_{0}(\bR^{d})$ and let $\zeta$ have
support in a ball $B$ of radius $R$ with center 
in $\bar G$ and satisfy $ 0 \leq\zeta\leq1$, 
  $|D\zeta|\leq K_{0}R^{-1}$,
$|D^{2}\zeta|\leq K_{0}^{2}R^{-2}$, where $K_{0}\geq1$ is a constant. Introduce
 $$
L'u:=Lu+cu,\quad
M u:=u  L'  \zeta +2a^{ij}D_{i}\zeta D_{j}u\quad (=L(\zeta u)-\zeta
Lu).
 $$
Then there exists a constant $N$, depending only
on $R_{0}$, $d,\delta,p,\|b\|$,  and $G$ (but not on $b_{0}$, $c_{0}$,
$r_{0}$, $\|c\|$, or
$\theta_{0}$), such that
$$
\|Mu\|_{L_{t}(G)}
\leq\varepsilon R^{-2\tau_{2}}\|D^{2}u\|_{L_{p}( G\cap B )}
$$
\begin{equation}
                                               \label{2.25.1}
+N(\varepsilon^{-\alpha}
K_{0}^{2}+\varepsilon^{-\beta}K_{0}^{2\gamma})R^{-2-2\tau_{2}}
\| u\|_{L_{p}( G\cap B )},
\end{equation}
where
$$
\alpha=\tau_{1}/(1-\tau_{1}),\quad \beta=\tau_{2}/(1-\tau_{2}),
\quad \gamma=(1-\tau_{2})^{-1} 
$$
and $ \tau_{1},\tau_{2}$ are specified in the proof.
\end{lemma}

Proof. Make the change of coordinates $y=x/R$ and, accordingly, set
  $u(x)=v(y)$. Under this change $B$ will be transformed into a ball
of radius one, the domain $G$ will also change, but, what is
important (due to $R\leq R_{0}$), the embedding theorems we need
in the transformed domain $B\cap G$ will hold with constants comparable to the
ones in the original $B\cap G$. Also observe that, if $Mu(x)=f(x)$, then
$$
v(y)\check L\check\zeta(y)+2\check a^{ij}(y)
D_{i}\check\zeta D_{j}v(y)=R^{2}f(Ry),
$$
where $\check L=\check a^{ij}(Ry)D_{ij}+Rb^{i}(Ry)D_{i}$,
$\check a^{ij}(y)=a^{ij}(Ry)$, $\check\zeta(y)=\zeta(Ry)$.
It is easy to check that the $L_{d}$-norm of the new $b$ remains the same.
It follows that we may concentrate on $R=1$.

In that case use H\"older's inequality
and embedding theorems (see, in particular,
Corollary 1.4.7/2 in \cite{Ma_85}). Observe that, 
$$
I :=\|u|b|\,|D\zeta|\,\|_{L_{t}( G )}\leq K_{0}\|b\|_{L_{d}( G )}
\|u\|_{L_{td/(d-t)}(
G\cap B )},
$$
and
since
$$
\frac{d}{p}-2<\frac{d(d-t)}{td} ,
$$
we have 
$$
I \leq \varepsilon\|D^{2}u\|_{L_{p}( G\cap B
)}+N\varepsilon^{-\tau_{1}/(1-\tau_{1})}K_{0}^{2}
\| u\|_{L_{p}( G\cap B )},
$$
where
$$
\tau_{1}= \frac{1}{2}\Big(1+ \frac{d}{p}-\frac{d}{t}\Big) .
$$

Also
$$
\frac{d}{p}-2<\frac{d}{t}-1,
$$
so that
$$
\||Du|\,|D\zeta|\,\|_{L_{t}( G )}\leq K_{0}\|Du\|_{L_{t}( G\cap B )}
$$
$$
\leq \varepsilon\|D^{2}u\|_{L_{p}( G\cap B )}+N
\varepsilon^{-\tau_{1}/(1-\tau_{1})}K_{0}^{2}
\| u\|_{L_{p}( G\cap B )}.
$$

Finally,
$$
\|| u|\,|D^{2}\zeta|\,\|_{L_{t}( G )}\leq K^{2}_{0}\|u\|_{L_{t}( G\cap B)} 
$$
$$
\leq \varepsilon\|D^{2}u\|_{L_{p}( G\cap B )}+N
\varepsilon^{-\tau_{2}/(1-\tau_{2})}K_{0}^{2/(1-\tau_{2})}
\| u\|_{L_{p}( G\cap B )},
$$
where
$$
\tau_{2}=(1/2)\Big(\frac{d}{p}-\frac{d}{t}\Big).
$$
Upon combining these estimates and observing that
$K_{0}\geq1$ and $\varepsilon\leq 1$,  we come to \eqref{2.25.1}
with $R=1$. The lemma is proved.

The following theorem allows us, in particular, to obtain
interior estimates.

\begin{theorem}
                                            \label{theorem 7.14.2} 
Under Assumptions \ref{assumption 3.1.1} and \ref{assumption 2.20.1}  ($p,r_{0}$)
let  
 $0<R\leq \text{\rm diam} (G)$, $z\in\bar{ G}$. Denote
$$
 G_{r}= G\cap B_{r}(z).
$$
Suppose that 
\begin{equation}
                                                       \label{10.24.1}
\zeta u\in \WO^{2}_{p}( G_{3R} )\quad\forall\zeta\in C^{\infty}_{0}
(B_{3R}(z) ),\quad
 Lu\in L_{p}( G_{3R}).
 \end{equation} 

Then
there exists a  constant  $N$, depending only on
 $d,\delta,p,  r_{0}$,  $b_{0}$, $c_{0}$, $G$, and $\|b\|$,   such that 
\begin{equation}
                                             \label{7.14.5}
\|u\|_{W^{ 2}_{p}( G_{R})}\leq
N \|Lu\|_{L_{p}( G_{2R})}+ N R^{-2}\|u\|_{ L_{p}( G_{2R})}.
\end{equation} 

\end{theorem}

Proof.  We may and will assume that $z=0$. In that case  set
$$
R_{m}=R\sum_{j=0}^{m}2^{-j},\quad D_{m}= G_{R_{m}},\quad m=0,1,2,....
$$
We need some functions
$\zeta_{m}\in C^{\infty}_{0}(\mathbb{R}^{d})$ such that
$\zeta_{m}(x)=1$ in $ B_{R_{m}}$, $\zeta_{m}(x)=0$ outside  
$ B_{R_{m+1}}$ and
$$
| D \zeta_{m}| \leq NR^{-1}2^{m} ,\quad
| D^{2}\zeta_{m}| \leq  NR^{-2}2^{ 2m},
 $$
where  $N=N(d)$.
To  construct them,
 take an infinitely differentiable function
 $h(t)$, $t\in(-\infty,\infty)$, such that
$h(t)=1$ for $t\leq0$, $h(t)=0$ for $t\geq1$ and
$0\leq h\leq1$. After this define
 $$
\zeta_{m}(x)=
h(2^{m+1}R^{-1}(|x|-R_{m})).
 $$

Now we put $u\zeta_{m}$ in \eqref{10.22.2} to
get
$$
\|u \|_{W^{ 2}_{p}(D_{m}) }\leq
\|u\zeta_{m}\|_{W^{ 2}_{p}( G)}
\leq N (\|L(u\zeta_{m})\|_{L_{p}( G)}
+\|u\zeta_{m}\|_{ L_{p}( G)})
$$
$$
\leq N\|Lu\|_{L_{p}( G_{2R})}
+\|M_{m}u\|_{L_{p}(G)}
+N\|u\|_{ L_{p}( G_{2R})}),
$$
where
 $$
M_{m} u:=u L'  \zeta_{m}+2a^{rs}D_{r}\zeta_{m}D_{s}u.
 $$
By Lemma \ref{lemma 2.25.1} with $t=p$ when $\tau_{2}=0$
(and $K_{0}\sim 2^{ m}$)
$$
\|M_{m}u\|_{L_{p}(G)}
\leq  (1/8)  \|u\zeta_{m+1}\| _{W^{ 2}_{p}( G)}
+N R^{-2}2^{ 2m}\|u\| _{ L_{p}( G_{2R})}.
$$

Then
$$
\|u\zeta_{m}\|_{W^{ 2}_{p}( G)}\leq 
N\|Lu\|_{L_{p}( G_{2R})}+
(1/8) \|u\zeta_{m+1}\| _{W^{ 2}_{p}( G)}
$$
$$
+ N  R^{-2}2^{2m}\|u\| _{ L_{p}( G_{2R})},
$$
$$
(1/8)^{m}\|u\zeta_{m}\|_{W^{ 2}_{p}( G)}\leq 
N(1/8)^{m}\|Lu\|_{L_{p}( G_{2R})}
$$
$$
+
(1/8)^{m+1}\|u\zeta_{m+1}\| _{W^{ 2}_{p}( G)}
+ N  R^{-2}2^{-m}\|u\| _{ L_{p}( G_{2R})}.
$$
By summing up over $m=0,1,..$ and cancelling like terms we obtain
$$
\|u\zeta_{m}\|_{W^{ 2}_{p}( G)}\leq 
 \|Lu\|_{L_{p}( G_{2R})}+ N  R^{-2} \|u\| _{ L_{p}( G_{2R})}.
$$
This proves \eqref{7.14.5} and the theorem.

In the following theorem we show that in our estimates
on the right one can have the $L_{p}$ norm of $u$ with lower $p$
(see \eqref{11.3.7}).

\begin{theorem}
                                                \label{theorem 11.3.4}
(i) Let    $q>d/2$, $q>p$,
 $0<R\leq \text{\rm diam} G $, $z\in\bar{ G}$.
  Denote
$$
 G_{r}= G\cap B_{r}(z).
$$

(ii) 
Introduce  
\begin{equation}
                                                         \label{3.11.6}
 \gamma =1+ \frac{2q-d}{d}\cdot\frac{p}{q}.
\end{equation}
Observe that $\gamma>1$ and introduce $p(n)=p\gamma^{n}$, $n=0,...,m-1$,
where $m-1$ is the largest $n$ such that $p(n)\leq q$.
Then set $p(m)=q$
and suppose that
  Assumption  \ref{assumption 3.1.1} is satisfied and Assumptions
\ref{assumption 2.20.1} ($ p(n),r_{0}$) are satisfied with the above $p(n)$'s,
$n=0,...,m$.

Then
\begin{multline}
                                                       \label{11.3.6}
\zeta u\in \WO^{2}_{p}( G_{2R} )\quad\forall\zeta\in C^{\infty}_{0}
(B_{2R}(z) ),\quad
 Lu\in  L_{q}( G_{2R})
\\ 
\Longrightarrow
\zeta u\in \WO^{ 2}_{q}( G_{2R} )\quad\forall\zeta\in C^{\infty}_{0}
(B_{2R}(z)).
\end{multline} 

Furthermore,
there exists a  constant  $N$, depending only
on $R$, $d,\delta,p,q$, $r_{0}$, $b_{0}$, $c_{0}$, 
 $\|b\|$,    and $G$,
such that, if the condition of the implication \eqref{11.3.6} holds, then
\begin{equation}
                                             \label{11.3.7}
\|u\|_{W^{2}_{q}( G_{R})}\leq
N(\|Lu\|_{ L_{q}( G_{2R})}+\|u\|_{ L_{p}( G_{2R})}).
\end{equation} 
\end{theorem}

Proof.  
Take $\lambda$ so large (see Theorem~\ref{theorem 10.17.1})
that $\lambda-L$
is invertible as an operator acting
from $\WO^{2}_{p(n)}( G)$ onto $ L_{p(n)}( G)$
for $n=0,...,m$.
 Also take a $u$ such that
the condition of the implication \eqref{11.3.6} holds, take 
a $\zeta\in C^{\infty}_{0}(B_{2R}(z))$, notice that
$\zeta u\in\WO^{2}_{p}( G)$ and denote  
$$
f=Lu,\quad
g=( L-\lambda)(\zeta u)=\zeta f+2a^{ij}u_{x^{i}}\zeta_{x^{j}}
+u(  L'   -\lambda)\zeta.
 $$

Observe that  for
$1\leq n<m$
\begin{equation}
                                                          \label{3.11.7}
\frac{d}{p(n-1)}-\frac{d}{p(n )}=
\frac{d}{p\gamma^{n}}(\gamma-1)\leq\frac{d}{p}(\gamma-1)
\leq \frac{2q-d}{q}=2-\frac{ d}{q}<1
\end{equation}
and $p  (n-1) \leq p(n )\leq q<d$. Therefore condition
\eqref{3.9.3} is satisfied with $t=p(n )$ and
$p(n -1 )$ in place of $p$. If $n=m$ and $p(m-1)=q$, then the
left-hand side of
\eqref{3.11.7} vanishes for $n=m$, and if $p(m-1)<q$, then
$q\leq p(m-1)\gamma$ and
$$
\frac{d}{p(m-1)}-\frac{d}{p(m )}=\frac{d}{p(m-1)}-\frac{d}{q}
\leq \frac{d}{p(m-1)}-\frac{d}{p(m -1)\gamma}<1.
$$
Hence, Lemma \ref{lemma 2.25.1} is applicable with
$t=p(n )$ and
$p(n -1 )$ in place of $p$ for any $n=1,...,m$.

It follows that  $g\in L_{p(n)}(G)$ if $\zeta u\in W^{2}_{p(n-1)}
 (G) $,
$n=1,...,m$.
By the choice of $\lambda$
the equation 
$$
( L-\lambda)w=g
$$
 has a solution in
$\WO^{2}_{p(1)}( G)
\subset\WO^{2}_{p}( G)$
which in addition is unique in $\WO^{2}_{p}( G)$.
Hence for $n=1$
\begin{equation}
                                               \label{11_14.5}
w=\zeta u\in\WO^{2}_{p(n)}( G)\quad\forall\zeta
\in C_{0}^{\infty}(B_{2R}).
\end{equation}

If $p(1)<q$, then
by repeating this argument with $p(1)$ in place of $p$,
we get \eqref{11_14.5} for $n=2$.
In this way we get this inclusion for all $n$ and 
this proves~\eqref{11.3.6}.

To prove \eqref{11.3.7}, we accompany the above argument
with estimates. By the choice of $\lambda$  and Lemma \ref{lemma 2.25.1},
for $n\geq1$ and any $\zeta,\eta\in C^{\infty}_{0}(B_{2R}(z))$
such that $\eta=1$ on the support of $\zeta$,
we have
$$
 \|\zeta u\|_{W^{2}_{p(n)}( G)}\leq
N\|\zeta f+2a^{ij}u_{x^{i}}\zeta_{x^{j}}
+u( L'   -\lambda)\zeta\|_{ L_{p(n)}( G)}
$$
$$
\leq N(\|f\|_{ L_{q}( G_{2R}) }+\|\eta u\|_{W^{2}_{p(n-1)}( G)}).
$$
By iterating the inequality between the extreme terms,
we obviously get that for any $\zeta\in C^{\infty}_{0}(B_{3R/2}(z))$
there is an $\eta\in C^{\infty}_{0}(B_{7R/4}(z))$ such that
$$
\|\zeta u\|_{W^{2}_{q}( G)}
\leq N(\|f\|_{ L_{q}( G_{2R} )}+\|\eta u\|_{W^{2}_{p}( G)}).
$$
Finally, recall that by Theorem \ref{theorem 7.14.2}
$$
\|\eta u\|_{W^{2}_{p}( G)}\leq N\|u\|_{W^{2}_{p}( G_{7R/4} )}
\leq N(\|f\|_{ L_{p}( G_{7R/2})}+\|u\|_{ L_{p}( G_{7R/2} )}).
$$ 
This yields \eqref{11.3.7} with $7R/2$ in place of $2R$.
However, obviously Theorem~\ref{theorem 7.14.2} is also true
with any number $>1$ in place of $2$. Then 
on the right 
in the above inequality
one can take $2R$ ($>7R/4$) in place of $7R/2$ and get \eqref{11.3.7}
in its original form. The theorem is proved.

\begin{remark}
                                            \label{remark 4.7.1}
The author does not know that, if $q=d/2$ and $p<q$,  the
implication
\eqref{11.3.6} is true or not. In that case a simple modification
of the above argument shows that $\zeta u\in \WO^{ 2}_{r}( G_{2R} )$
with any $r<d/2$. The reason behind this restriction
lies in Theorem \ref{theorem 10.17.1} in which the case $d/2=q=p$
is not allowed.
\end{remark}

\begin{corollary}
                                                \label{corollary 3.5.1}
Under the assumptions of Theorem \ref{theorem 11.3.4}
if $u\in W^{2}_{p,\loc}(G)$ satisfies $Lu=0$ in $G$,
then $u\in  W^{2}_{q,\loc}(G)$. In particular, $u\in C^{2-d/q}_{\loc}
(G)$.
\end{corollary}

Below we use the constant $d_{0}=d_{0}(d,\delta, \|b\| )\in(d/2,d)$
introduced in  \cite{Kr_19_1}.
From   Corollary 6.3 of \cite{Kr_20}
and Corollary \ref{corollary 3.5.1} we obtain the following.
\begin{corollary}[Harnack inequality]
                                                \label{corollary 3.7.1}
Under the assumptions of Theorem \ref{theorem 11.3.4},
let $q\geq d_{0}$, $R\in(0,\infty]$ and let $u\in W^{2}_{p}(B_{2R})$ be a
nonnegative   function satisfying $Lu=0$ (a.e.)
in $B_{2R}$ with $c\equiv0$. Then for any $x,y\in B_{R}$ we have   
$u(x)\leq Nu(y)$, where $N=N(d,\delta,\|b\|)$.
\end{corollary}

It would be interesting to know if this result can be obtained
by using purely PDE methods as in \cite{Sa_10}.

The following theorem will be used when $c\equiv0$,
so that we can take $q$ as close to $d$ as we like.

\begin{theorem}
                                              \label{theorem 12.14.1}
 Under Assumption \ref{assumption 3.1.1} 
suppose that either  (a) $q>d/2$ and the condition (ii) 
 of Theorem \ref{theorem 11.3.4}
is satisfied or (b) $p>d/2$.  
Set $\bar \lambda= \lambda_{1}(d,\delta,p,r_{0} , b_{0} , c_{0},   G  )$
where the latter is introduced in Theorem \ref{theorem 10.17.1}.

 Then 
 there exists an integer $m_{0}$,
depending only on $p$  and $d$, and  a constant  $N$,
depending only on $d,\delta,p,q ,r_{0} $, $b_{0}, c_{0}, 
\|b\|,  G$, such that for any $f\in  L_{p}( G)$ we have
 \begin{equation}
                                                       \label{12.13.6}
\sup_{x\in G}| R_{\bar\lambda+c} ^{m_{0}}f(x)|\leq
N\|f\|_{ L_{p}( G)}.
\end{equation}  
\end{theorem}

Proof. If $p>d/2$, we are done due to Theorem \ref{theorem 10.17.1}
and embedding theorems.
In case $p\leq d/2$ we also have $p<q$ and we use
$p(n)=\gamma^{n}p$, $n=0,1,...,m-1$, $p(m)=q$, and set
$$
u_{0}=f,\quad
 u_{n}= R_{\bar\lambda+c} ^{n}f,\quad n\geq1.
$$
Observe  that, for $n\geq0$, we have
$$
\bar\lambda  u_{n+1}-Lu_{n+1}= u_{n},
$$
so that $u_{n+1}\in\WO^{2}_{p}
( G)$ and
$$
\|u_{n+1}\|_{W^{2}_{p}( G)}
\leq N\|u_{n}\|_{ L_{p}( G)}
\leq N\|u_{n}\|_{ L_{p(n)}( G)}.
$$
By Theorem \ref{theorem 11.3.4}
 $$
\|u_{n+1}\|_{W^{2}_{p(n)}( G)}
\leq N\|u_{n}\|_{ L_{p(n)}( G)}+N
\|u_{n+1}\|_{ L_{p }( G)}.
$$
Hence
\begin{equation}
                                         \label{07.12.27.3}
\|u_{n+1}\|_{W^{2}_{p(n)}( G)}
\leq N\|u_{n}\|_{ L_{p(n)}( G)},
\end{equation} 
and by embedding theorems  
$$
\|u_{n+1}\|_{ L_{p(n+1)}( G)}
\leq N\|u_{n}\|_{ L_{p(n)}( G)}.
$$
Iterating this yields  that for $n\geq0$
\begin{equation}
                                         \label{2.25.5}
\|u_{n}\|_{ L_{p(n)}( G)}
\leq N\|u_{0}\|_{ L_{p(0)}( G)}= 
N\|f\|_{ L_{p}( G)},
\end{equation}
where the constants 
$N$ depend on the data as in the statement of the theorem
and they also depend on $n$. Now we fix an $n=n(p,d)$ ($\leq m$) so that $p(n)
>d/2$ and from \eqref{07.12.27.3} and \eqref{2.25.5} and embedding theorems
conclude that 
$$
\sup_{x\in G}|u_{n+1}(x)|\leq N\|u_{n+1}
\|_{W^{2}_{p(n)}( G)}
\leq N\|f\|_{ L_{p}( G)},
$$
which shows that \eqref{12.13.6} holds with $m_{0}=n+1$.
 The theorem is proved.

\mysection{Two auxiliary results using probability theory}

Here we assume that the coefficients of $L$
 satisfy only Assumption \ref{assumption 3.1.1}
and  are infinitely differentiable.
   Reading this section requires
some acquaintance with the basic notions of
stochastic calculus which are found, for instance, in \cite{Kr_02}.
For the reader's orientation we sketch some of them. A $d$-dimensional
Wiener precess $w_{t}$ is the mathematical model of Brownian motion
and is a continuous random process with independent increments and
independent coordinates such that $w^{i}_{t}-w^{i}_{s}$
has normal distribution with zero mean and variance $|t-s|$
for any $i=1,...,d$. It\^o proved that one can define the stochastic
integral
$$
\int_{0}^{t}f_{t}\,dw_{t}
$$
for random $\bR^{d}$-valued $f_{t}$
as the limit of usual integral sums provided 
that
$f$, say,  is measurable bounded and, for each $t$, $f_{t}$ and 
the process $w_{s+t}-w_{t}$,
$s>0$, are independent. After that, by using Perron's method of successive
approximations, he
showed that under our above assumptions on $a$ and $b$, for any $x$, the equation
\begin{equation}
                                    \label{2.29.4}
x_{t}=x+\int_{0}^{t}\sqrt{2a(x_{s})}\,dw_{s}+\int_{0}^{t}
b(x_{s})\,ds
\end{equation}
has a unique solution such that  for each $t$, $x_{t}$ and 
the process $w_{s+t}-w_{t}$,
$s>0$, are independent. Finally, what we need is It\^o's formula,
which implies (see \cite{Kr_77}) that if $D$ is a bounded domain
 $u\in
W^{2}_{d}(D)\cap C(\bar D)$ and $c_{t}\geq 0$   is measurable bounded and,
for each
$t$,
$c_{t}$ and  the process $w_{s+t}-w_{t}$,
$s>0$, are independent, then for any $x\in D$
$$
u(x)=Ee^{-\phi_{\tau}}u(x_{\tau})+E\int_{0}^{\tau}e^{-\phi_{t}}
(c_{t}u(x_{t})-Lu(x_{t}))\,dt,
$$
where $x_{t}$ is the solution of \eqref{2.29.4}, $\tau$
is its first exit time from $D$ and
$$
\phi_{t}=\int_{0}^{t}c_{s}\,ds.
$$

\begin{lemma}
                                                 \label{lemma 2.23.1}
Let $\lambda\geq\nu>0$ and let $u\in \WO^{2}_{d}(G)$
satisfy $ \lambda u-Lu \leq1$ in $G$. Then $\lambda u \leq \mu$,
where $\mu<1$ is a constant depending only on $\nu, d,\delta$, $\|b\|$,
and the diameter of $G$.

\end{lemma}

Proof.   By It\^o's formula 
$$
u(x)=E\int_{0}^{\tau}f(x_{t})e^{-\lambda t}\,dt,
$$
where $\tau$ is the first exit time of $x_{t}$ from $G$
and $f=\lambda u-Lu$. It follows that
$$
\lambda u(x) \leq E\int_{0}^{\tau}\lambda e^{-\lambda t}\,dt
=1-Ee^{-\lambda \tau}\leq 1-Ee^{-\nu \tau}.
$$
By Corollary 2.7 of \cite{Kr_19} there exist constants
$\kappa=\kappa(d,\delta,\|b\|,\text{diam}(G))>0$
and $N=N(d,\delta,\|b\|)$ such that, for any $T>0$
\begin{equation}
                                    \label{2.23.3}
P(\tau>T)\leq Ne^{-\kappa T}.
\end{equation}
Hence,
$$
Ee^{-\nu \tau}\geq e^{-\nu T}P(\tau\leq T)\geq e^{-\nu T}(1-
Ne^{-\kappa T})
$$
and $\lambda u(x) \leq 1-  e^{-\nu T}(1-
Ne^{-\kappa T})=:\mu$, where $\mu<1$ for an appropriate choice of $T$.
The lemma is proved.

\begin{lemma}
                                                 \label{lemma 2.29.7}
Let $\lambda\geq\nu>0$ and let $u\in W^{2}_{d}(B_{R})$
satisfy $ \lambda u-Lu \leq 0$ in $B_{R}$. Then 
\begin{equation}
                                                  \label{2.29.8}
u(0)\leq 2e^{-\kappa\sqrt{\nu} R}\max_{\partial B_{R}}u_{+} ,
\end{equation}
where $\kappa=\kappa(d,\delta,\|b\|)>0$.

\end{lemma}

This lemma is a direct corollary of Theorem 2.10 of \cite{Kr_19} because
in light of It\^o's formula
$$
u(0)= E\Big(e^{-\lambda
\tau-\phi_{\tau}}u (x_{\tau})+\int_{0}^{\tau}e^{-\lambda t-\phi_{t}} (\lambda
-L)u(x_{t})\,dt\Big)\leq 
 E e^{-\lambda \tau}\max_{\partial B_{R}}u_{+},
$$
where 
$$
\phi_{t}=\int_{0}^{t}c(x_{s})\,ds,
$$
$\tau$ is the first exit time of $x_{t}$ from $B_{R}$
and $x_{t}$ is the solution of \eqref{2.29.4} with $x=0$.

 \mysection{Solvability of $\lambda u-Lu=f$ in $G$ for $\lambda\geq 0$}
 
Here we follow the line of arguments 
from Section 11.3 of \cite{Kr_08}. We take $p,q$   
as in Section \ref{section 3.11.1}.
We also take any $q'>d/2,q'>p$, produce $\gamma $
by using \eqref{3.11.6} with $q'$ in place of $q$,
then introduce $m$ and $p(n)$ as in the statement of Theorem
\ref{theorem 11.3.4}   with $q'$ in place of $q$.
We also take $\lambda\geq 0$.  

\begin{assumption}
                                          \label{assumption 3.12.2}
 Assumption  \ref{assumption 3.1.1} is satisfied and
either (a) $q=p$ (which is only possible if $p>d/2$) and 
Assumption \ref{assumption
2.20.1} ($ p ,r_{0}$) is satisfied or (b) $q>p$ and
Assumptions \ref{assumption 2.20.1} ($ p(n),r_{0}$)
  are satisfied with the above $p(0),...,p(m)$.
\end{assumption}

This assumption is supposed to be satisfied throughout   
the section.
 
\begin{theorem}
                                                 \label{theorem 12.13.4}
There exists
a constant $N$ depending only on $ d,\delta $,
$p$, $q$, $r_{0}$, $b_{0}$, $c_{0}$,   $\|b\|$, and $G$,
such that for any   $u\in\WO^{2}_{p}( G)$
\begin{equation}
                                                       \label{12.14.2}
\|u\|_{W^{2}_{p}( G)}\leq N\|(\lambda-L)u\|_{ L_{p}( G)}.
\end{equation}
Furthermore, for any $f\in L_{p}(G)$ there exists a unique
$u\in\WO^{2}_{p}( G)$ such that $\lambda u-Lu=f$ in $G$.
\end{theorem}

Proof. In light of the method of continuity it suffices to prove
the first assertion, while proving which we may
assume that the coefficients of $L$ are infinitely
differentiable.
If $\lambda\geq \bar\lambda$,
with $\bar\lambda$ taken from Theorem \ref{theorem 12.14.1},  
 the result
is known from   Theorem \ref{theorem 10.17.1}. Therefore
we will only concentrate on  
 $
0\leq\lambda<\bar\lambda$. 
Define   
$$
f=\lambda u-Lu
$$
so that
$$
\bar\lambda u-Lu=(\bar\lambda-\lambda)u+f,\quad
u=(\bar\lambda-\lambda)R_{\bar\lambda+c}u+R_{\bar\lambda+c}f,
$$
and by induction on $n$
$$
u=[(\bar\lambda-\lambda)R_{\bar\lambda+c}]^{n }u
+\sum_{i=0}^{n-1}[(\bar\lambda-\lambda)R_{\bar\lambda+c}]^{i}
R_{\bar\lambda+c}f,
$$
where $n$ is any integer $\geq1$. We thus have the beginning
of the   Neumann series.

Introduce the constants $N_{1}$ and $M_{n}$
  so that
$$
  \|R_{\bar\lambda}g\|_{ L_{p}( G)}\leq N_{1}
\|g\|_{ L_{p}( G)}\quad\forall
g\in  L_{p}( G),\quad
M_{n}=\sum_{i=0}^{n-1}  \bar\lambda ^{i}N_{1}^{i+1}.
$$
Finally, let $| G|$ be the volume of $ G$ and
take $m_{0}$
from Theorem \ref{theorem 12.14.1} when $c\equiv0$ and $q$ is replaced  
by $q'=p(m)$.
For $n>m_{0}$, in light of  Corollary \ref{corollary 3.10.1}  
$$
\|u\|_{ L_{p}( G)}
 \leq| G|^{1/p} \bar\lambda^{n}\sup_{x\in G}
 R_{\bar\lambda }^{n-m_{0}}R_{\bar\lambda }^{ m_{0}}|u|(x) 
+M_{n}\|f\|_{ L_{p}( G)}.
$$
By Lemma \ref{lemma 2.23.1} the above supremum is dominated by
$$
\bar\lambda^{m_{0}-n}\mu^{n-m_{0}}\sup_{x\in G} R_{\bar\lambda }^{ m_{0}}|u|(x) ,
$$
where $\mu<1$, which by Theorem \ref{theorem 12.14.1}
is less than
$$
N_{2}\bar\lambda^{m_{0}-n}\mu^{n-m_{0}}\|u\|_{ L_{p}( G)}.
$$
Hence,
$$
\|u\|_{ L_{p}( G)}
 \leq N_{2}| G|^{1/p} \bar\lambda ^{m_{0}}
\mu^{n-m_{0}}\|u\|_{ L_{p}( G)}+M_{n}\|f\|_{ L_{p}( G)}.
$$
We fix $n$ so that  $N_{2}| G|^{1/p} \bar\lambda ^{m_{0}}
\mu^{n-m_{0}}\leq 1/2$ and then arrive at 
$$
\|u\|_{ L_{p}( G)}\leq 2M_{n}\|f\|_{ L_{p}( G)}.
$$
Now  to get
\eqref{12.14.2} it only remains to refer to Remark \ref{remark 2.21.1}.
The theorem is proved.

\begin{corollary}[Maximum principle]
                                              \label{corollary 11.7.1}
Let   $u\in \WO^{2}_{p}(G)$.
Then
\begin{equation}
                                                   \label{11.7.7}
\|u_{\pm}\|_{ L_{p}( G)}\leq N\|(\lambda u-Lu)_{\pm}\|_{ L_{p}( G)}
\end{equation} 
where $N$ depends only on $ \delta,d$,
$p$, $q$, $r_{0}$, $b_{0}$, $c_{0}$, $\|b\|$, and the diameter of $G$.
In particular, if $u\in \WO^{2}_{p}(G)$
and $Lu-\lambda u\geq0$ in $G$, then $u\leq0$
in $G$.
\end{corollary}

This corollary is derived from Theorem \ref{theorem 12.13.4}
in the same way as Theorem 11.3.3 of \cite{Kr_08}
is derived from Theorem 11.3.2.

\mysection{Equations in the whole space with $\lambda$ large}

  We take $p,q$   
as in Section \ref{section 3.11.1}.
The following is a slight restatement of part of Theorem 6.4.1
of \cite{Kr_08}.
\begin{lemma}
                          \label{lemma 2.19.20}
Let $s\in (1,\infty)$.
If Assumptions \ref{assumption 3.1.1}
 and \ref{assumption 2.20.1} ($ s ,r_{0}$)  are  satisfied,
then there exist
  $\lambda_{0}=\lambda_{0}(d,\delta,s)$
and   $N_{0}=N_{0}(d,\delta,s,r_{0})$ such that, for 
any $u\in W^{2}_{s}(\bR^{d})$
and $\lambda\geq \lambda_{0}$,
\begin{equation}
                                 \label{2.19.30s}
\|D^{2}u\|_{L_{s}(\bR^{d})}+\sqrt\lambda
\|D u\|_{L_{s}(\bR^{d})}+ \lambda\| u\|_{L_{s}(\bR^{d})}  
\leq N_{0} \|L_{0}u-\lambda u\|_{L_{s}(\bR^{d})}.
\end{equation}
 Furthermore, for any $f\in L_{s}(\bR^{d})$
there exists a unique $u\in W^{2}_{s}(\bR^{d})$
such that $L_{0}u-\lambda u=f$.
\end{lemma}

In this section we suppose that 
 Assumptions \ref{assumption 3.1.1} and \ref{assumption 2.20.1}
($ p ,r_{0}$) are satisfied.
Below by $\lambda_{0}$ we mean the one from Lemma
\ref{lemma 2.19.20} for $s=p$.

\begin{theorem}  
                                                       \label{theorem 2.20.10}
Introduce $N^{*} =N^{*}(p,d)$ as the best
constant such that
$$
\|Du\|_{L_{pd/(d-p)}(\bR^{d})}\leq N^{*}\|D^{2}u\|_{L_{p}(\bR^{d})}. 
$$
for any $u\in W^{2}_{p}(\bR^{d})$. 
Assume that
\begin{equation}     
                            \label{2.20.30}
2N_{0}N^{*}(\|b\|_{L_{d}(\bR^{d})}+\|c\|_{L_{q}(\bR^{d})})\leq1 .
\end{equation}

Then for any $u\in W^{2}_{p}(\bR^{d})$
and $\lambda\geq \lambda_{0}$,
\begin{equation}
                                 \label{2.20.40}
\|D^{2}u\|_{L_{p}(\bR^{d})}+\sqrt\lambda
\|D u\|_{L_{p}(\bR^{d})}+\lambda\| u\|_{L_{p}(\bR^{d})}
\leq 2N_{0} \|L u-\lambda u\|_{L_{p}(\bR^{d})}.
\end{equation}
 Furthermore, for any $f\in L_{p}(\bR^{d})$
there exists a unique $u\in W^{2}_{p}(\bR^{d})$
such that $L u-\lambda u=f$.
\end{theorem}

The proof of this theorem is achieved by repeating that
of Theorem \ref{theorem 2.20.1}.

\begin{remark}
                                                  \label{remark 2.20.20}
Similarly to Remark \ref{remark 2.20.2} we observe that
 the operator $L$ is bounded
as an operator from $W^{2}_{p}(\bR^{d})$ to $L_{p}(\bR^{d})$
as long as $b\in L_{d}(\bR^{d})$ and $c\in L_{q}(\bR^{d})$.
\end{remark}  
 
Next, for our fixed $b$ and $c$ there exists a $b_{0},c_{0}\geq 0$ such that  
\begin{equation}     
                            \label{2.20.50}
2N_{0}N^{*}(\|bI_{|b|\geq b_{0}}\|_{L_{d}(\bR^{d})}
+\|cI_{c\geq c_{0}}\|_{L_{q}( \bR^{d})})\leq1 .
\end{equation}
 
Obviously we may take the same $b_{0}$, $c_{0}$
in \eqref{2.20.5} and \eqref{2.20.50}.

\begin{theorem}
                                     \label{theorem 3.5.10}
There exist
  $\lambda_{1}\geq 1, N $,
depending only on $d,\delta,p,r_{0},b_{0}$, and $c_{0}$,    
  such that, for any $u\in W^{2}_{p}(\bR^{d})$
and $\lambda\geq \lambda_{1}$,
\begin{equation}
                                 \label{3.1.10}
\|D^{2}u\|_{L_{p}(\bR^{d})}+\sqrt\lambda\| Du\|_{L_{p}(\bR^{d})}+\lambda\|
u\|_{L_{p}(\bR^{d})}
\leq N  \|L u-\lambda u\|_{L_{p}(\bR^{d})}.
\end{equation}
 Furthermore, for any $f\in L_{p}(\bR^{d})$
there exists a unique $u\in W^{2}_{p}(\bR^{d})$
such that $L u-\lambda u=f$.
\end{theorem}

One proves this theorem in the same way as Theorem \ref{theorem 10.17.1}.

 We denote the solution from Theorem \ref{theorem 3.5.10}
by $R_{\lambda+c}f$.
 
\begin{remark}
                                                  \label{remark 2.21.10}

By taking   $\lambda=\lambda_{1}$ in \eqref{3.1.10}
we see that for the same kind of $N$ as in \eqref{3.1.10}
and any $u\in W^{2}_{p}( \bR^{d})$
\begin{equation}
                                               \label{10.22.20}
\| u\|_{W^{2}_{p}(\bR^{d})}\leq N\big(\|Lu\|_{L_{p}( \bR^{d} )}
+\|u\|_{L_{p}( \bR^{d})} \big).
\end{equation}
\end{remark}

\mysection{Equations in the whole space with $\lambda$ small}  

 We take $p,q$   
as in Section \ref{section 3.11.1}.
We also take any $q'>d/2,q'>p$, produce $\gamma $
by using \eqref{3.11.6} with $q'$ in place of $q$,
then introduce $m$ and $p(n)$ as in the statement of Theorem
\ref{theorem 11.3.4}   with $q'$ in place of $q$.
We also take $\lambda > 0$.

\begin{assumption}
                                          \label{assumption 3.12.1}
 Assumption     
\ref{assumption 3.1.1} is satisfied and
either (a) $q=p$ and Assumption \ref{assumption 2.20.1} ($ p ,r_{0}$)
is satisfied or (b) $q>p$ and
Assumptions  \ref{assumption 2.20.1} ($ p(n),r_{0}$)
  are satisfied with the above $p(0),...,p(m)$.
\end{assumption}

 Recall that $B_{R}=\{x:|x|<R\}$. Take the constant $\kappa$
from Lemma \ref{lemma 2.29.7} and define $ R'
=R'(d,\delta,\|b\|,\lambda) \geq 4$ so that
$$
2e^{-\kappa\sqrt\lambda( R'-2)  }\leq 1/2.
$$
\begin{lemma}
                                               \label{lemma 11.3.1}
Under Assumption \ref{assumption 3.12.1}
let $u$ and $f $ be bounded infinitely differentiable
functions. Assume that $f=0$  outside
 $B_{1}$ and $\lambda u-Lu=f$ in $\bR^{d}$. Also assume that the
coefficients of $L$ are infinitely differentiable.  

 Then there exists a  constant  
$N$, depending only on $\lambda$, 
 $d,\delta,p,q, r_{0}$, $b_{0}$, $c_{0}$,   and $\|b\|$, such
that
$$
\|u/v\|_{ L_{p}(\bR^{d} )}\leq N\|f\|_{ L_{p}(\bR^{d} )},
$$
where $v(x)=e^{-\kappa\sqrt\lambda|x|}$.
\end{lemma}
 
Proof. We follow the proof of Lemma 11.6.1 of \cite{Kr_08}.
 Relying on classical
results, define $h\in W^{2}_{q}(B_{R'})$  as a unique solution
of 
$$
\lambda h-Lh=0\quad\text{in}\quad B_{R'}  \quad\text{such that}
\quad  w:=h-u\in
\WO^{2}_{q}(B_{R'}).
$$
 By regularity results $h$
is infinitely differentiable in $\bar{B}_{R'}$
and $h=u$ on $\partial B_{R'}$. Hence
 $w$ is infinitely differentiable in $\bar{B}_{R'}$,
vanishes on $\partial B_{R'}$,
 and satisfies 
$$
\lambda w-Lw= f.
$$ 

Notice that $\lambda u-Lu=0$ outside  $B_{1}$ and by the maximum principle
$$
|u(x)|\leq  
\max_{|x|=2}|u|\quad\text{for}
\quad|x|\geq 2.
$$

Taking this into account, taking $x$ as the new origin, and using
 Lemma \ref{lemma 2.29.7},  
 we obtain
\begin{equation}
                                                 \label{11.3.3}
|u(x)|\leq 2e^{-\kappa\sqrt\lambda(|x|-2)}
\max_{|x|=2}|u|\quad\text{for}
\quad|x|\geq 2.
\end{equation} 
Also observe that 
by the maximum principle
$$
|h|\leq \max_{|x|=R'}|u|
$$
 in $B_{R'}$.

Now we claim that to prove the lemma, it suffices to prove that
\begin{equation}
                                                 \label{11.3.5}
|w(x)|\leq N\| f\|_{ L_{p}(B_{R'} )}\quad\text{for}
\quad|x|=2.
\end{equation} 

Indeed, if  \eqref{11.3.5} holds, then
$$
\max_{|x|=2}|u| \leq
\max_{|x|=2}| h |+\max_{|x|=2}|w|\leq \max_{|x|=R'}|u|
+ N\|f\|_{ L_{p}(\bR^{d} )}
$$
$$
\leq 2e^{-\kappa\sqrt\lambda(R'-2)  } \max_{|x|=2}|u|
+ N\|f\|_{ L_{p}(\bR^{d} )},
$$
which for our choice of $R' $ yields
$$
\max_{|x|=2}|u|\leq N\|f\|_{ L_{p}(\bR^{d} )}.
$$

Coming back to \eqref{11.3.3}  
and using that $e^{2\kappa\sqrt \lambda}\leq N$
we get that
$$
 \|u/v\|_{ L_{p}(B_{2}^{c})}\leq N\|f\|_{ L_{p}(\bR^{d} )}.
$$
The remaining part of the norm is also bounded by
$N\|f\|_{ L_{p}(\bR^{d} )}$ since $|u|\leq|h|+|w|$,
$$
\max_{B_{R'}}|h| \leq\max_{|x|=R'}|u|
\leq   \max_{|x|=2}|u|
\leq N\|f\|_{ L_{p}(\bR^{d} )} ,
$$
   and by Theorem \ref{theorem 12.13.4} we have
$$
\|w\|_{ L_{p}(B_{ R' })}\leq N\|f\|_{ L_{p}(B_{ R' })}.
$$ 
Thus, indeed we need only prove \eqref{11.3.5}.

By the maximum principle $|w|\leq \psi$, where $\psi$
is a $\WO^{2}_{q}(B_{R'})$-solution of $L\psi=-|f|$. 
So it suffices to estimate $\psi$ on $|x|=2$.
  Take a point $x_{0}$ with $|x_{0}|=2$
and observe that by embedding theorems
we have
$$
|\psi(x_{0})|\leq N\|\psi\|_{W^{2}_{q}(B_{1/2}(x_{0}))}.
$$
Next, we  
use the local regularity result
from Theorem \ref{theorem 11.3.4}. Then we find
$$
\|\psi\|_{W^{2}_{q}(B_{1/2}(x_{0}))}
\leq N \| L\psi\|_{ L_{q}(B_{1}(x_{0}))}+N\|\psi\|
_{ L_{p}(B_{1}(x_{0}))} .
$$
Here the first term on the right is zero since $ f=0$
  outside of $B_{1}$ and the second term
is less than $N\|f\|_{ L_{p}(B_{R'})}$ by Theorem \ref{theorem
12.13.4}. The lemma is proved.

The above proof of Lemma \ref{lemma 11.3.1}
is slightly different from the proof of  Lemma 11.6.1 of \cite{Kr_08} 
and is drift-specific
because we needed to use Lemma \ref{lemma 2.29.7}, whose counterpart
in \cite{Kr_08} was obtained by using simple barriers.
Contrary to that the following theorem is derived
from  Lemma \ref{lemma 11.3.1} by literally repeating
the derivation  of Theorem 11.6..2 of \cite{Kr_08}
from Lemma 11.6.1 of \cite{Kr_08}.

\begin{theorem}
                                                \label{theorem 11.5.2}
Under Assumption \ref{assumption 3.12.1}
for any $f\in  L_{p}(\bR^{d})$
there exists a unique $u\in W^{2}_{p}(\bR^{d})$ such that
$\lambda u-Lu=f$. Moreover, there exists a constant $N$, depending only
$\lambda$, 
 $d,\delta,p, q, r_{0}$, $b_{0}$, $c_{0}$, and $\|b\|$, such that
$$
\|u\|_{ W^{2}_{p}(\bR^{d} )}\leq N\|f\|_{ L_{p}(\bR^{d} )}.
$$
\end{theorem}

One more result concerning elliptic equations
which will be proved in the next section is the following
stability theorem in which Assumption 
\ref{assumption 3.12.1} is not imposed.

\begin{theorem}
                                                \label{theorem 3.8.1}
 Let $q=p\geq d_{0}$, where $d_{0}=d_{0}(d,\delta, \|b\| )\in(d/2,d)$
is taken from    \cite{Kr_19_1}, and suppose that Assumptions
 \ref{assumption 3.1.1} and \ref{assumption 2.20.1} ($ p ,r_{0}$) are
satisfied.
Let  $a^{n},b^{n},c^{n}$, $n= 1,2,...$,
be   sequences of smooth bounded functions with values in the
set of symmetric $d\times d$ matrices having all
 eigenvalues in $[\delta, \delta^{-1}]$, in
$\bR^{d}$, and in $[0,\infty)$, respectively,
 such that 
$ a^{n}\to a$ on $\bR^{d}$ (a.e.) and
$$
\|b-b^{n}\|_{L_{d}(\bR^{d})}+\|c^{n}-c\|_{L_{q}(\bR^{d})}\to0
$$
 as $n\to\infty$.
Take $\lambda>0$, $f\in L_{q}(\bR^{d})$, and introduce $u^{n}$
as  unique  $W^{2}_{q}(\bR^{d})$-solutions of $\lambda
u^{n}-L^{n}u^{n}=f
$, where the operators $L^{n}$ are constructed from  $a^{n},b^{n},c^{n}$.
Then at each point of $\bR^{d}$ we have $u^{n}\to u$
as $n\to\infty$, where $u\in W^{2}_{q}(\bR^{d})$
is a unique solution of $\lambda u -L  u =f  $.
\end{theorem}

The author does not know if this theorem holds for $q\in ( d/2,d_{0})$.
In this range we have  unique solvability in $W^{2}_{p}(\bR^{d})$
with $1<p\leq q$, but could it happen that there is no stability?

\mysection{Weak uniqueness of solutions of stochastic equations}

Here  we let $q=p\geq d_{0}$, where
$d_{0}=d_{0}(d,\delta, \|b\| )\in(d/2,d)$ is taken from  
 \cite{Kr_19_1}, and suppose that Assumptions
 \ref{assumption 3.1.1} and \ref{assumption 2.20.1} ($ p ,r_{0}$) are
satisfied.
Take $x\in\bR^{d}$. Recall that according to Theorem 1.1 of \cite{Kr_19_1}
there exists 
a probability space $(\Omega ,\cF ,P )$,
a filtration of $\sigma$-fields $\cF _{t}\subset \cF $, $t\geq0$,
a process $w _{t}$, $t\geq0$, which is a $d$-dimensional Wiener process
relative to $\{\cF _{t}\}$, and an $\cF _{t}$-adapted
process $x_{t}$ such that 
 (a.s.) for all   $t\geq0$
\begin{equation}
                                                 \label{11.29.2}
x _{t}=x  +\int_{0}^{t}\sqrt{2a (x_{s})}\,dw_{s}
+\int_{0}^{t}b (x_{s}) \,ds.
\end{equation}

Take $f\in L_{q}(\bR^{d})$ and $\lambda>0$.
By Theorem \ref{theorem 11.5.2} there is  a unique 
$u\in W^{2}_{q}(\bR^{d})$ such that
$\lambda u-Lu=f$. By Theorem 1.3 of \cite{Kr_19_1}
It\^o's formula is applicable so that
for all $t\geq0$   and $L_{s}=a^{ij}(x_{s})D_{ij}+b^{i}(x_{s})
D_{i}$ and $\sigma_{s}=\sqrt{2a (x_{s})}$ 
\begin{equation}
                                              \label{9.6.6}
u(x_{t})=u(x)+\int_{0}^{t}L_{s}u(x_{s})\,ds
+\int_{0}^{t}D_{i}u(x_{s})\sigma^{ik}_{s}\,dw^{k}_{s}
\end{equation}
and the last term is a square integrable martingale.
By Theorem 1.5 of \cite{Kr_19_1}
$$
E\int_{0}^{\infty}e^{-\lambda t}c(x_{t})\,dt<\infty.
$$
Therefore, It\^o's formula is applicable to
\begin{equation}
                                                 \label{4.8.4}
u(x_{t})\exp(-\lambda t-\int_{0}^{t}c(x_{s})\,ds).
\end{equation}\begin{remark}
                                    \label{remark 4.8.1}
Applying It\^o's
formula to \eqref{4.8.4} yields that
$$
u (x)=E\int_{0}^{\infty}f(x _{t})\exp(-\lambda t
-\int_{0}^{t}c (x _{s})\,ds)\,dt.
$$
\end{remark}

These facts and the standard argument based on considering
resolvent operators (see, for instance, the arguments in 
\cite{Kr_19_1} after Theorem 1.1 there) immediately proves
the following weak uniqueness theorem.

\begin{theorem}
                                              \label{theorem 3.8.10}
All solutions of \eqref{11.29.2} on all possible
probability spaces have the same distribution
on $C([0,\infty),\bR^{d})$.
\end{theorem}

{\bf Proof of Theorem \ref{theorem 3.8.1}}. Let $x^{n}_{t}$
be solutions of
\begin{equation}
                                                 \label{4.8.3}
x^{n} _{t}=x  +\int_{0}^{t}\sqrt{2a^{n} (x^{n}_{s})}\,dw_{s}
+\int_{0}^{t}b^{n} (x^{n}_{s}) \,ds
\end{equation}
on the same probability space as $x_{t}$ or on different
ones.
By Theorem 1.1 of  \cite{Kr_19_1}
the set of distributions of $x^{k }_{\cdot}$ on 
$C([0,\infty),\bR^{d})$ is tight and any weakly converging
subsequence of distributions converges weakly to 
the distribution of one of solutions of
\eqref{11.29.2}, which is the only one in light
of Theorem  \ref{theorem 3.8.10}. Hence,
the whole sequence of distributions  
of $x^{n }_{\cdot}$ weakly converges to the distribution
of $x_{\cdot}$. In particular, for any  $n_{0}$ and smooth bounded $g$

$$
\lim_{n\to\infty}E\int_{0}^{\infty}g(x^{n}_{t})\exp(-\lambda t
-\int_{0}^{t}c^{n_{0}}(x^{n}_{s})\,ds)\,dt
$$
\begin{equation}
                                                 \label{4.8.5}
= E\int_{0}^{\infty}g(x _{t})\exp(-\lambda t
-\int_{0}^{t}c^{n_{0}}(x _{s})\,ds)\,dt.
\end{equation}

At this point it is appropriate to mention that by It\^o's
formula
$$
u^{n}(x)=E\int_{0}^{\infty}f(x^{n}_{t})\exp(-\lambda t
-\int_{0}^{t}c^{n}(x^{n}_{s})\,ds)\,dt.
$$

Next, since $\|b^{n}\|_{L_{d}(\bR^{d})}\leq  \|b \|_{L_{d}(\bR^{d})}+1$
for sufficiently large $n$, by Theorem 1.5 of \cite{Kr_19_1}
for any $\lambda>0$, $r\geq d_{0}$,
and   $g(x)$ given on $\bR^{d}$  we have
\begin{equation}
                                          \label{9.5.6}
E\int_{0}^{\infty}e^{-\lambda t} 
|g(x^{n }_{t})|\,dt\leq N \lambda^{d/(2r)-1}\|g\|_{L_{r}(\bR^{d})},
\end{equation}
where $N$ is independent of $f$ and $n$.
Below all constants like this one are called $N$.
The same estimate holds for $x_{t}$ in place of $x^{n }_{t}$.
Hence also taking into account
\eqref{4.8.5} we get that, for $\varepsilon>0$ and smooth bounded $g$
such that $\|f-g\|_{L_{q}(\bR^{d})}\leq\varepsilon$,
$$
\nliminf_{n\to\infty}u^{n}(x)
\geq -N\varepsilon
+\nliminf_{n\to\infty}E\int_{0}^{\infty}g(x^{n}_{t})\exp(-\lambda t
-\int_{0}^{t}c^{n}(x^{n}_{s})\,ds)\,dt
$$
$$
\geq -N\varepsilon+
E\int_{0}^{\infty}g(x _{t})\exp(-\lambda t
-\int_{0}^{t}c^{n_{0}}(x _{s})\,ds)\,dt
$$
$$
-\sup|g|E\int_{0}^{\infty}e^{-\lambda t}\Big(\int_{0}^{t}
|c^{n}(x^{n}_{s})-c^{n_{0}}(x^{n} _{s})|\,ds\Big)dt.
$$
Integrating by parts we see that the last expectation equals
$$
\lambda^{-1} E\int_{0}^{\infty}e^{-\lambda t} 
|c^{n}(x^{n}_{t})-c^{n_{0}}(x^{n} _{t})| \,dt\leq N\|c^{n}-c^{n_{0}}\|
_{L_{q}(\bR^{d})}.
$$
It follows that
$$
\nliminf_{n\to\infty}u^{n}(x)
\geq -N\varepsilon-N\sup|g|\|c -c^{n_{0}}\|
_{L_{q}(\bR^{d})}
$$
$$
+
E\int_{0}^{\infty}g(x _{t})\exp(-\lambda t
-\int_{0}^{t}c^{n_{0}}(x _{s})\,ds)\,dt.
$$
By using similar estimates for $u(x)$ we conclude that
$$
\nliminf_{n\to\infty}u^{n}(x)
\geq -N\varepsilon-N\sup|g|\|c -c^{n_{0}}\|
_{L_{q}(\bR^{d})}+u(x).
$$
By letting $\varepsilon\downarrow 0$ and $n_{0}\to\infty$
we arrive at
$$
\nliminf_{n\to\infty}u^{n}(x)\geq u(x).
$$
This result is also true if we replace  $f$ with $-f$
and this, certainly, proves the theorem.

{\bf Acknowledgment}. The author is thankful to E. Priola
for some comments and to Hongjie Dong who conjectured
that the results of the kind presented in the article
should be true  and kindly pointed out to the
author some glitches in the earlier version 
of the article and a simple error which led to estimating nonexisting terms
albeit not affecting the results. The author   also brings his thanks
to two referees whose comments helped improve and sometimes correct
the presentation.

\end{document}